\documentclass[11pt]{amsart2000}
\usepackage{amssymb}

\begin{document}
\setlength{\unitlength}{0.01in}
\linethickness{0.01in}
\begin{center}
\begin{picture}(474,66)(0,0)
\multiput(0,66)(1,0){40}{\line(0,-1){24}}
\multiput(43,65)(1,-1){24}{\line(0,-1){40}}
\multiput(1,39)(1,-1){40}{\line(1,0){24}}
\multiput(70,2)(1,1){24}{\line(0,1){40}}
\multiput(72,0)(1,1){24}{\line(1,0){40}}
\multiput(97,66)(1,0){40}{\line(0,-1){40}}
\put(143,66){\makebox(0,0)[tl]{\footnotesize Proceedings of the Ninth Prague Topological Symposium}}
\put(143,50){\makebox(0,0)[tl]{\footnotesize Contributed papers from the symposium held in}}
\put(143,34){\makebox(0,0)[tl]{\footnotesize Prague, Czech Republic, August 19--25, 2001}}
\end{picture}
\end{center}
\vspace{0.25in}
\setcounter{page}{191}
\title{A survey of $J$-spaces}
\author{E. Michael}
\address{University of Washington\\
Seattle, WA, U.S.A.}
\thanks{The author was an invited speaker at the Ninth Prague Topological 
Symposium.}
\thanks{This contribution is excerpted from a published article.
Reprinted from 
Topology and its Applications,
Volume 102, number 3,
E. Michael,
$J$-Spaces,
pp. 315--339,
Copyright (2000),
with permission from Elsevier Science \cite{3}.}
\thanks{E. Michael,
{\em A survey of $J$-spaces},
Proceedings of the Ninth Prague Topological Symposium, (Prague, 2001),
pp.~191--193, Topology Atlas, Toronto, 2002}
\begin{abstract}
This note is a survey of $J$-spaces.
\end{abstract}
\maketitle

\section{Basic concepts}

A space $X$ is a \emph{$J$-space} if, whenever $\{A,B\}$ is a closed cover
of $X$ with $A\cap B$ compact, then $A$ or $B$ is compact. 
A space $X$ is a \emph{strong $J$-space} if every compact $K\subset X$ is
contained in a compact $L\subset X$ with $X\backslash L$ connected.
[As in \cite{3}, all maps are continuous and all spaces are Hausdorff.]

\subsection{}
Every strong $J$-space $X$ is a $J$-space. 
The two concepts coincide when $X$ is locally connected, but in general
(even for closed subsets of ${\mathbb R}^2$) they do not.

\section{Examples}

\subsection{}
A topological linear space $X$ is a (strong) $J$-space if and only if
$X\ne{\mathbb R}$.

\subsection{}
If $X$ and $Y$ are connected and non-compact, then $X\times Y$ is a strong
$J$-space.\footnote{This was proved in \cite{4}.}

\subsection{}
Let $Y$ be a compact manifold with boundary $B$, and let $A\subset B$. 
Then $Y\backslash A$ is a (strong) $J$-space if and only if $A$ is connected.

\section{Characterizations by closed maps} 

A map $f:X\to Y$ is called \emph{boundary-perfect} if $f$ is closed and
$\operatorname{bdry} f^{-1}(y)$ is compact for every $y\in Y$. It follows
from \cite{2}
that every closed map $f:X\to Y$ from a paracompact space $X$ to a
$q$-space $Y$ is boundary-perfect.\footnote{$q$-spaces (see 
\cite{2}) include all locally compact and all metrizable spaces.}

\subsection{}
A space $X$ is a $J$-space if and only if every boundary-perfect map
$f:X\to Y$ onto a non-compact space $Y$ is perfect.

\subsection{}
If $X$ is a $J$-space, then every boundary-perfect map $f:X\to Y$ has at
most one non-compact fiber. 
The converse holds if $X$ is locally compact.

\subsection{}
Let $X$ be paracompact and locally compact. 
Then the following are equivalent. 
\begin{itemize}
\item[(a)]
$X$ is a $J$-space.
\item[(b)]
Every closed map $f:X\to Y$ onto a non-compact, locally compact space $Y$
is perfect.
\item[(c)]
Every closed map $f:X\to Y$ onto a locally compact space $Y$ has at most
one non-compact fiber.
\end{itemize}

\subsection{}
Let $X$ be metrizable. 
Then the following are equivalent
\begin{itemize}
\item[(a)]
$X$ is a $J$-space.
\item[(b)]
Every closed map $f:X\to Y$ onto a non-compact, metrizable space $Y$ is
perfect.
\end{itemize}

\section{Characterization by compactifications} 

Call a set $A\subset Y$ a \emph{boundary set} for $Y$ if 
$\operatorname{Int} A=\emptyset$ and, whenever $U\supset A$ is open in $Y$
and 
$\{W_1, W_2\}$ is a disjoint, relatively open cover of $U\backslash A$,
then no $y\in A$ lies in $\overline W_1\cap \overline W_2$. 
Call a set $A\subset Y$ a \emph{strong boundary set} for $Y$ if 
$\operatorname{Int} A=\emptyset$ and, whenever $U\supset A$ is open in
$Y$, then every
$y\in A$ has an open neighborhood $V\subset U$ with $V\backslash A$
connected.

It is easy to see that, if $Y$ is a manifold with boundary $B$, then
every $A\subset B$ is a strong boundary set for $Y$. And it follows from
the proof of \cite[Lemma 4]{1} (or from \cite[Proposition
3.5]{Michael}) that, if $Y$ is completely regular, then $\beta
X\backslash X$ is a boundary set for $\beta X$.

\subsection{}
Let $Y$ be a compactification of $X$, and suppose either that $X$ is
locally compact or that $Y$ is metrizable. 
Then the following are equivalent.
\begin{itemize}
\item[(a)]
$X$ is a (strong) $J$-space.
\item[(b)]
$Y\backslash X$ is connected and a (strong) boundary set for $Y$.
\end{itemize}

\section{Preservation}

\subsection{}
$J$-spaces are preserved by boundary-perfect images. 
(False for strong $J$-spaces, even with perfect images.)

\subsection{}
$J$-spaces and strong $J$-spaces are preserved by monotone, perfect
pre-images.

\subsection{}
If $X_1$, $X_2$ are connected, then $X_1\times X_2$ is a 
(strong) $J$-space if and only if either $X_1$, $X_2$ are both
(strong) $J$-spaces or both are non-compact. 

\subsection{}
Let $\{X_1, X_2\}$ be a closed cover of $X$ with $X_1\cap X_2$ compact. 
Then $X$ is a (strong) $J$-space if and only if $X_1$, $X_2$ are both
(strong) $J$-spaces and $X_1$ or $X_2$ is compact.

\subsection{}
If $X$ is a (strong) $J$-space, so is every component of $X$. 
(False for $J$-spaces).


\begin{thebibliography}{1}

\bibitem{1}
Melvin Henriksen and J.~R. Isbell, \emph{Local connectedness in the {S}tone-\v
  {C}ech compactification}, Illinois J. Math. \textbf{1} (1957), 574--582.
  \MR{20 \#2688}

\bibitem{Michael}
E.~Michael, \emph{Cuts}, Acta Math. \textbf{111} (1964), 1--36. \MR{29 \#5233}

\bibitem{2}
\bysame, \emph{A note on closed maps and compact sets}, Israel J. Math.
  \textbf{2} (1964), 173--176. \MR{31 \#1659}

\bibitem{3}
\bysame, \emph{${J}$-spaces}, Topology Appl. \textbf{102} (2000), no.~3,
  315--339. \MR{2002b:54019}

\bibitem{4}
Krzysztof Nowi{\'n}ski, \emph{Closed mappings and the {F}reudenthal
  compactification}, Fund. Math. \textbf{76} (1972), no.~1, 71--83. \MR{48
  \#2978}

\end{thebibliography}
\providecommand{\bysame}{\leavevmode\hbox to3em{\hrulefill}\thinspace}
\providecommand{\MR}{\relax\ifhmode\unskip\space\fi MR }
\providecommand{\MRhref}[2]{%
  \href{http://www.ams.org/mathscinet-getitem?mr=#1}{#2}
}
\providecommand{\href}[2]{#2}

\end{document}